\documentclass[10pt]{article}
\usepackage{amsmath,amsthm,amsfonts,amssymb,latexsym,graphicx}

\usepackage{xcolor}

\usepackage{algorithm}
\usepackage[noend]{algpseudocode}
\algrenewcommand\algorithmicthen{\relax}
\algdef{C}[IF]{IF}{ElsIf}[1]{\textbf{elif}\ #1}
\algrenewcommand\algorithmicdo{\relax}

\usepackage[colorlinks=true,citecolor=blue,pdfpagemode=UseNone,pdfstartview=FitH]{hyperref}

\allowdisplaybreaks[4]

\newcommand{\st}{\mid}

\renewcommand{\d}{\,\mathrm{d}}

\newcommand{\E}{\mathbb{E}}
\renewcommand{\P}{\mathbb{P}}
\newcommand{\R}{\mathbb{R}}
\newcommand{\N}{\mathbb{N}}

\newcommand{\FDP}{\textup{FDP}}

\newcommand{\AAA}{\mathcal{A}}
\newcommand{\FFF}{\mathcal{F}}

\theoremstyle{plain}
\newtheorem{theorem}{Theorem}

\newtheorem{lemma}[theorem]{Lemma}
\newtheorem{proposition}[theorem]{Proposition}

\theoremstyle{definition}

\theoremstyle{remark}
\newtheorem{remark}[theorem]{Remark}

\title{Tuning the false discovery rate in the closed eBH procedure}
\author{Vladimir Vovk}
\date{November 21, 2026}

\begin{document}
\maketitle

  \begin{abstract}
    The closed eBH procedure controls the false discovery rate
    under arbitrary dependence between the input e-values.
    It takes the target false discovery rate $\alpha$ as its input
    and outputs a rejection set $R$,
    the usual goal being to maximize its size $\left|R\right|$.
    This note discusses a modification of the closed eBH procedure
    that allows choosing the target false discovery rate post hoc:
    the rejection set and its level $\alpha$ can be chosen after seeing the data,
    and the expected ratio of the false discovery proportion to $\alpha$ is at most~$1$.
    We adapt Vovk and Wang's (2023) algorithm for computing discovery e-matrices
    to computing, in quadratic time,
    the frontier of achievable combinations
    $(\left|R\right|,\alpha)$ of the size of a rejection set $R$ and its level $\alpha$.
    The frontier allows us to compute an optimal rejection set
    in view of a given utility function $F$,
    $F(\left|R\right|,\alpha)\to\max$,
    striking a balance between the goals of reducing $\alpha$ and increasing $\left|R\right|$.

    \medskip

    \noindent
    The version of this note at \url{http://alrw.net/e} (Working Paper 11)
    is updated most often.
  \end{abstract}

\section{Introduction}

The closed eBH procedure in \cite[(7)]{Xu/etal:arXiv2509}
(see also \cite[Definition~9.31]{Ramdas/Wang:2025})
is an elegant improvement of the original eBH procedure \cite{Wang/Ramdas:2022},
but its disadvantage is that we need to decide in advance
on the target false discovery rate $\alpha$.
For example, what if we can slightly exceed $\alpha$ but then reject many more null hypotheses?
This note discusses a computationally efficient way of choosing both the rejection set $R$
and its \emph{level} $\alpha$ that plays the role of the target false discovery rate.

We start in Sect.~\ref{sec:closed} with fixing our terminology and notation
and a short description of the closed eBH procedure
and its standard property of validity.
The next section, Sect.~\ref{sec:relative},
gives a straightforward post hoc modification of the closed eBH procedure,
which we call the ``relative eBH procedure'';
it is used, more or less explicitly, in several publications
(see, e.g., \cite[Sect.~8.3]{Xu/etal:arXiv2509} and \cite[Sect.~9.6]{Ramdas/Wang:2025}).
The levels $\alpha$ output by the relative eBH procedure are valid in relative terms:
the expected ratio of the false discovery proportion to the reported level is at most 1.
Section~\ref{sec:algorithm} gives an efficient algorithm for the relative eBH procedure
by adapting an algorithm in \cite[Algorithm~4]{Vovk/Wang:2023}
for computing discovery e-matrices.
For both closed and relative
eBH procedures we give large-deviation inequalities for sequences of studies,
and Sect.~\ref{sec:conclusion} summarizes some advantages of the relative procedure.

We will introduce most of our terminology and notation
at the point where we need them.
The notation for the natural numbers is $\N:=\{1,2,\dots\}$ and $\N_0:=\{0,1,\dots\}$.
The terms ``increasing'' and ``decreasing'' will be used in the wide sense
allowing intervals of constancy.

\section{The closed eBH procedure}
\label{sec:closed}

Our sample space $(\Omega,\AAA)$ is arbitrary,
and our statistical model is the set $\mathfrak{P}(\Omega)$
of all probability measures on $(\Omega,\AAA)$;
since this model contains all probability measures on $\Omega$,
there is no real loss of generality.
For $\theta\in\mathfrak{P}(\Omega)$ and a nonnegative random variable $E$,
we set $\E^{\theta}(E):=\int E\d\theta$
(this is always defined, perhaps with $\E^{\theta}(E)=\infty$).

Suppose that we are given $K$ \emph{e-variables} $E^{(1)},\dots,E^{(K)}$
for testing hypotheses $H_1,\dots,H_K\subseteq\mathfrak{P}(\Omega)$,
which are our \emph{base hypotheses}:
each $E^{(k)}$ is a nonnegative random variable satisfying
$\E^{\theta}(E^{(k)})\le1$ for all $\theta\in H_k$.
We would like to reject some of the base hypotheses
(in fact, as many of them as possible under a validity constraint).
A \emph{rejection set} is a set $R\subseteq[K]$, where $[K]:=\{1,\dots,K\}$;
this is the set of base hypotheses, represented by their indices,
that the researcher chooses to reject.
For each $\theta\in\mathfrak{P}(\Omega)$, let
\[
  I_{\theta}
  :=
  \{k\in[K]\st\theta\in H_k\}
\]
be the indices of the hypotheses containing $\theta$.
If the researcher rejects $H_k$, we refer to this decision as a \emph{discovery}.
We say that the discovery is \emph{true} if $\theta\notin H_k$,
and it is \emph{false} if $\theta\in H_k$,
where $\theta$ is the true (unknown) probability measure governing the data generation.

Let
\begin{equation*}
  \FDP_I(R)
  :=
  \frac{\left|R\cap I\right|}{\left|R\right| \vee 1},
\end{equation*}
where $I,R\subseteq[K]$.
In particular, $\FDP_{I_{\theta}}(R)$ is the \emph{false discovery proportion} of a rejection set $R$
when the data are generated from $\theta$.
Our informal goal is to make $\FDP_{I_{\theta}}(R)$ small for a large $R$;
the obvious obstacle is, of course, that we do not know $\theta$.

For a nonempty $I\subseteq[K]$, the arithmetic mean
\begin{equation}\label{eq:E_I}
  E_I
  :=
  \frac{1}{\left|I\right|}
  \sum_{k\in I}
  E^{(k)}
\end{equation}
is an e-variable for testing $\cap_{k\in I}H_k$.
Following \cite[Sect.~9.6]{Ramdas/Wang:2025}, we also set $E_{\emptyset}:=1$.
Therefore, $\E^{\theta}(E_{I_{\theta}})\le1$ for all $\theta\in\mathfrak{P}(\Omega)$.
We will sometimes refer to $(E_I)_{I\subseteq[K]}$ as our \emph{e-collection}.

The \emph{closed eBH procedure} of \cite[(7)]{Xu/etal:arXiv2509}
(see also \cite[Definition 9.31]{Ramdas/Wang:2025})
takes $\alpha\in(0,1]$ (the target false discovery rate) as an input
and outputs a (largest) set $R$ such that
\begin{equation}\label{eq:closed}
  \alpha E_I\ge\FDP_I(R)
  \quad
  \text{for all $I\subseteq[K]$}.
\end{equation}

The validity of the closed eBH procedure follows immediately from $\E^{\theta}(E_{I_{\theta}})\le1$.
The rejection set in the following proposition does not have to be a largest one:
the researcher may choose any rejection set satisfying \eqref{eq:closed},
possibly after looking at the data.

\begin{proposition}\label{prop:FDR}
  Let $\alpha\in(0,1]$.
  For every $\theta\in\mathfrak{P}(\Omega)$
  and every random set $\hat R\subseteq[K]$ satisfying \eqref{eq:closed},
  \[
    \E^{\theta}
    \left(
      \FDP_{I_{\theta}}(\hat R)
    \right)
    \le
    \alpha.
  \]
\end{proposition}

\begin{proof}
  Taking $I:=I_{\theta}$ in \eqref{eq:closed} gives $\FDP_{I_{\theta}}(\hat R)\le\alpha E_{I_{\theta}}$,
  and we know that $\E^{\theta}(E_{I_{\theta}})\le1$.
\end{proof}

The left-hand side of the inequality in Proposition~\ref{prop:FDR}
is the false discovery rate of $\hat R$ under $\theta$,
and so Proposition~\ref{prop:FDR} says that the closed eBH procedure
controls the false discovery rate at level $\alpha$ \cite{Xu/etal:arXiv2509}.

Proposition~\ref{prop:FDR} is about a single study.
Let us now consider a sequence of studies indexed by $n=1,2,\dots$,
each of the kind described above
(with the number of base hypotheses that may vary from study to study).
In the $n$th study, the researcher chooses a rejection set $R_n$,
possibly after looking at the data,
and reports a level $\alpha_n$ for it.
The studies are \emph{sequential}:
the base hypotheses and the e-variables of the $n$th study,
as well as the probability measure $\theta_n$ generating its data,
may depend on the outcomes of the earlier studies.
Let $\FFF_n$ be the $\sigma$-algebra generated by the outcomes of the first $n$ studies,
including the researcher's choices
($\FFF_0$ is the trivial $\sigma$-algebra),
and let $\P$ and $\E$ refer to the joint distribution of all studies.
We assume that $\theta_n$ is the conditional distribution of the data of the $n$th study given $\FFF_{n-1}$.
Let $I_n:=I_{\theta_n}$ be the set of indices of the true base hypotheses in the $n$th study,
$E_n$ be the e-variable $E_{I_n}$ of the $n$th study
(i.e., \eqref{eq:E_I} for $I:=I_n$),
and $\FDP_n:=\FDP_{I_n}(R_n)$ be its false discovery proportion.
Our assumptions imply
\begin{equation}\label{eq:sequential}
  \E(E_n\mid\FFF_{n-1})\le1,
  \quad
  n\in\N,
\end{equation}
and this is all we will need;
\eqref{eq:sequential} holds, in particular, when the studies are independent.

\begin{proposition}\label{prop:closed-sequential}
  Let $\alpha_{\max}\in(0,1]$.
  Suppose that the level $\alpha_n\in(0,\alpha_{\max}]$ for the $n$th study
  is chosen before the $n$th study is carried out
  (formally, it is $\FFF_{n-1}$-measurable)
  and that $R_n$ satisfies \eqref{eq:closed} with $\alpha_n$ in place of $\alpha$
  (and with the e-variables of the $n$th study).
  Then, for all $N\in\N$ and $\epsilon>0$,
  \begin{equation}\label{eq:LD-fixed}
    \P
    \left(
      \frac1N
      \sum_{n=1}^N
      (\FDP_n-\alpha_n)
      \ge
      \epsilon
    \right)
    \le
    \exp
    \left(
      -\frac{N\epsilon^2}{2\alpha_{\max}+2\epsilon/3}
    \right),
  \end{equation}
  and
  \[
    \limsup_{N\to\infty}
    \frac1N
    \sum_{n=1}^N
    (\FDP_n-\alpha_n)
    \le
    0
    \quad
    \text{a.s.}
  \]
\end{proposition}

\begin{proof}
  Taking $I:=I_n$ in \eqref{eq:closed} gives $\FDP_n\le\alpha_n E_n$,
  i.e., $\FDP_n-\alpha_n\le\alpha_n(E_n-1)$;
  besides, $\FDP_n-\alpha_n\le1$ since $\FDP_n\le1$.
  It remains to apply Lemma~\ref{lem:Bernstein} in Appendix~\ref{app:Bernstein}
  (with $Y_n:=\FDP_n-\alpha_n$, $c_n:=\alpha_n$, and $c:=\alpha_{\max}$).
\end{proof}

In particular, if the closed eBH procedure is used at the same level $\alpha$ in all studies,
the average false discovery proportion of the first $N$ studies
is at least $\alpha+\epsilon$
with probability at most $\exp(-N\epsilon^2/(2\alpha+2\epsilon/3))$,
and its $\limsup$ as $N\to\infty$ is at most $\alpha$ a.s.
These are frequentist interpretations of the closed eBH procedure's validity.
The rejection sets $R_n$ may be chosen after looking at the data,
but only among the sets satisfying \eqref{eq:closed}
at the level $\alpha_n$ chosen before the study;
the next section relaxes this requirement.

See Algorithm~\ref{alg:frontier}
for an efficient implementation of the closed eBH procedure.

\section{The relative eBH procedure}
\label{sec:relative}

In this section we allow the level $\alpha$ to be chosen post hoc.
For each rejection set $R\subseteq[K]$
define the random variable
\begin{equation}\label{eq:alpha}
  \alpha_R
  :=
  \min\{\alpha\in[0,1]\st\forall I\subseteq[K]:\alpha E_I\ge\FDP_I(R)\},
\end{equation}
where $\min\emptyset:=\infty$.
We will refer to the version of the closed eBH procedure
that assigns the level $\alpha_R$ to each rejection set $R$
as the \emph{relative eBH procedure}.

Suppose we have decided on a utility function $F:\N_0\times[0,1]\to\R$
that is increasing in its first argument and decreasing in the second,
and extend it by $F(r,\infty):=-\infty$, $r\in\N_0$.
We will assume that it is computable in time $O(1)$.
The main topic of this note is to design an efficient algorithm for the optimization problem
\begin{equation}\label{eq:objective}
  F(\left|R\right|,\alpha_R)\to\max
\end{equation}
over $R\subseteq[K]$.
The algorithm is given as its input the e-values $E^{(k)}(\omega)$, $k\in[K]$,
for the realized outcome $\omega\in\Omega$.

By the definition \eqref{eq:alpha} of $\alpha_R$,
\begin{equation}\label{eq:basic}
  \FDP_I(R)/\alpha_R
  \le
  E_I
\end{equation}
for all $I,R\subseteq[K]$, with the conventions $0/0:=0$ and $c/\infty:=0$
(if $\alpha_R=0$, then $\FDP_I(R)=0$ for all $I$).
The following property of validity is stated in \cite[Theorem~9.32]{Ramdas/Wang:2025}
and \cite[Sect.~8.3]{Xu/etal:arXiv2509}.

\begin{proposition}\label{prop:validity}
  For every $\theta\in\mathfrak{P}(\Omega)$ and every random set $\hat R\subseteq[K]$,
  in particular for a solution $\hat R$ of the optimization problem \eqref{eq:objective},
  \[
    \E^{\theta}
    \left(
      \frac{\FDP_{I_{\theta}}(\hat R)}{\alpha_{\hat R}}
    \right)
    \le
    \E^{\theta}
    \left(
      \max_{R\subseteq[K]}
      \frac{\FDP_{I_{\theta}}(R)}{\alpha_R}
    \right)
    \le
    1.
  \]
\end{proposition}

\begin{proof}
  The inequality \eqref{eq:basic} with $I:=I_{\theta}$ gives
  $\max_R\FDP_{I_{\theta}}(R)/\alpha_R\le E_{I_{\theta}}$,
  and we know that $\E^{\theta}(E_{I_{\theta}})\le1$.
\end{proof}

Proposition~\ref{prop:validity} says that $\alpha_{\hat R}$ is a valid post hoc level for $\hat R$
in the sense of \cite{Grunwald:2024} and \cite{Koning:arXiv2312}:
the expected ratio of the false discovery proportion to the reported level is at most~$1$,
even though $\hat R$, and therefore its level, is chosen after looking at the data.
Since $\hat R$ can be any random set,
the utility function $F$ in \eqref{eq:objective} does not even have to be fixed in advance:
one can compute the whole frontier (defined formally in the next section)
and choose a point on it after inspecting it.

\begin{remark}
  Proposition~\ref{prop:validity} does not imply that the expected false discovery proportion
  is at most the expected reported level,
  $\E^{\theta}(\FDP_{I_{\theta}}(\hat R)-\alpha_{\hat R})\le0$.
  For example, suppose that only the first of $K:=2$ base hypotheses is true
  and that the e-variables $(E^{(1)},E^{(2)})$ take the values $(2,100)$ and $(0,100)$
  with probability $1/2$ each.
  If the researcher chooses the largest rejection set whose level is at most $1/4$
  and reports its level,
  she reports $\{1,2\}$ with level $1/4$ in the first case
  and $\{2\}$ with level $1/50$ in the second case.
  The expected difference between the false discovery proportion and the reported level is then
  \[
    \frac{1/2-1/4}{2} + \frac{0-1/50}{2}
    =
    \frac{23}{200}
    >
    0,
  \]
  whereas the expected ratio is $(2+0)/2=1$.
  (For checking calculations such as the preceding one
  and the one in Remark~\ref{rem:alpha-star} below,
  it is useful to consult Lemma~\ref{lem:frontier}.)
  The most obvious way to ensure that the expected false discovery proportion
  is at most the expected reported level
  is to fix, before seeing the data, a reference level $c\in(0,1]$
  and to set
  \begin{equation*}
    \alpha^c_R
    :=
    \max_{I\subseteq[K]}
    \bigl(
      \FDP_I(R)-c(E_I-1)
    \bigr);
  \end{equation*}
  this \emph{absolute eBH procedure}
  will then satisfy $\E^{\theta}(\FDP_{I_{\theta}}(\hat R)-\alpha^c_{\hat R})\le0$.
  While the absolute eBH procedure contains
  the closed eBH procedure with the target false discovery rate of $c$
  as an admissible option (similarly to the relative eBH procedure),
  it might not contain much more in practically interesting cases
  as $\alpha^c_R\ge(1-E^{\min})c$ for any rejection set $R$,
  where $E^{\min}$ is the smallest e-value among $E^{(1)},\dots,E^{(K)}$.
\end{remark}

Proposition~\ref{prop:validity} is about a single study.
Let us now consider a sequence of studies, as in Sect.~\ref{sec:closed},
where in each study the researcher chooses both the rejection set $R_n$
and a level $\alpha_n$ for it, allowing $\alpha_n>\alpha_{R_n}$.
We have the following analogue of Proposition~\ref{prop:closed-sequential}
for levels chosen post hoc.

\begin{proposition}\label{prop:sequential}
  Suppose that $\alpha_{R_n}\le\alpha_n\le1$ for all $n\in\N$,
  where $\alpha_{R_n}$ is the level \eqref{eq:alpha} of $R_n$ in the $n$th study,
  and that, for some $\alpha_{\min}\in(0,1]$,
  $\alpha_n\ge\alpha_{\min}$ whenever $R_n\ne\emptyset$.
  Then, for all $N\in\N$ and $\epsilon>0$,
  \begin{equation}\label{eq:LD-relative}
    \P
    \left(
      \frac1N
      \sum_{n=1}^N
      \frac{\FDP_n}{\alpha_n}
      \ge
      1+\epsilon
    \right)
    \le
    \exp
    \left(
      -\frac{\alpha_{\min}N\epsilon^2}{2+2\epsilon/3}
    \right),
  \end{equation}
  and
  \[
    \limsup_{N\to\infty}
    \frac1N
    \sum_{n=1}^N
    \frac{\FDP_n}{\alpha_n}
    \le
    1
    \quad
    \text{a.s.}
  \]
  (with $0/0:=0$).
\end{proposition}

A natural choice for $\alpha_n$ in Proposition~\ref{prop:sequential}
is $\alpha_n:=\max(\alpha_{R_n},\alpha_{\min})$;
this ensures the existence of a positive floor for $\alpha_n$
required by statements of the form~\eqref{eq:LD-relative}
while using the optimal level $\alpha_{R_n}$ above the floor.

\begin{proof}[Proof of Proposition~\ref{prop:sequential}]
  Set $Y_n:=\alpha_{\min}(\FDP_n/\alpha_n-1)$.
  If $R_n=\emptyset$, then $Y_n=-\alpha_{\min}$.
  If $R_n\ne\emptyset$, then $\FDP_n/\alpha_n\le\FDP_n/\alpha_{R_n}\le E_n$ by \eqref{eq:basic},
  and $\FDP_n/\alpha_n\le1/\alpha_{\min}$.
  In both cases,
  \[
    Y_n
    \le
    \alpha_{\min}
    \left(
      \min(E_n,1/\alpha_{\min})
      -
      1
    \right)
    \le
    \min(1,\alpha_{\min}(E_n-1)),
  \]
  and it remains to apply Lemma~\ref{lem:Bernstein}
  with $c_n:=c:=\alpha_{\min}$
  and with $\alpha_{\min}\epsilon$ in place of $\epsilon$.
\end{proof}

The rejection sets $R_n$ and the levels $\alpha_n$ in Proposition~\ref{prop:sequential}
can be chosen in any way, in particular post hoc
(as long as $\alpha_{R_n}\le\alpha_n$),
and so Proposition~\ref{prop:sequential} gives an interpretation of post hoc levels
in terms of the ratios $\FDP_n/\alpha_n$,
provided the reported levels are bounded away from $0$.
For a level $\alpha$ fixed in advance ($\alpha_n:=\alpha$ for all $n$),
we can take $\alpha_{\min}:=\alpha$,
and then \eqref{eq:LD-relative} with $\epsilon/\alpha$ in place of $\epsilon$
becomes \eqref{eq:LD-fixed} with $\alpha_n=\alpha$.

\section{Computing the frontier}
\label{sec:algorithm}

In
this section we fix the realized outcome $\omega\in\Omega$
and write $e_k:=E^{(k)}(\omega)$, $e_I:=E_I(\omega)$,
and $\alpha_R:=\alpha_R(\omega)$.
Set
\[
  \beta_R
  :=
  \max_{I\subseteq[K]}
  \frac{\FDP_I(R)}{e_I},
  \quad
  R\subseteq[K],
\]
with the conventions $0/0:=0$ and $c/0:=\infty$ for $c>0$.
Then
\[
  \alpha_R
  =
  \begin{cases}
    \beta_R & \text{when $\beta_R\le1$}\\
    \infty & \text{otherwise}.
  \end{cases}
\]
The \emph{frontier} is the function $r\mapsto\alpha^*_r$, where
\[
  \alpha^*_r
  :=
  \min_{R\subseteq[K]:\left|R\right|=r}
  \alpha_R,
  \quad
  r\in\{0,\dots,K\}.
\]
Since $F$ is decreasing in its second argument,
\begin{equation}\label{eq:reduction}
  F(\left|R\right|,\alpha_R)
  \le
  F(\left|R\right|,\alpha^*_{\left|R\right|})
  \quad
  \text{for all $R\subseteq[K]$}.
\end{equation}
Therefore, to solve the optimization problem \eqref{eq:objective}
it suffices to know the frontier and, for each $r$,
a set of size $r$ whose level is $\alpha^*_r$.

\begin{remark}\label{rem:alpha-star}
  The sequence $\alpha^*_r$ is not necessarily monotone:
  e.g., for the e-values $(1,10,10)$ we have
  $(\alpha^*_0,\dots,\alpha^*_3)=(0,2/11,1/7,1/3)$,
  and so rejecting two base hypotheses requires a smaller level than rejecting one.
  For solving optimization problems \eqref{eq:objective}
  we only need, of course, the Pareto frontier,
  i.e., the running minimum of our frontier from the right.
\end{remark}

If we do not know anything about the nature of the base hypotheses,
it makes sense to reject a number of them with the largest e-values
(see, e.g., \cite[Sect.~6]{Vovk/Wang:2023} and \cite[Lemma~52]{Xu/etal:arXiv2509}),
and Lemma~\ref{lem:frontier} below shows that this is optimal in our current setting.
From now on
we assume, without loss of generality,
that the e-values are sorted in ascending order,
\begin{equation}\label{eq:sorted}
  e_1\le\dots\le e_K,
\end{equation}
and set $R_0:=\emptyset$ and
\[
  R_r
  :=
  \{K-r+1,\dots,K\},
  \quad
  r\in[K];
\]
$R_r$ consists of the hypotheses with the $r$ largest e-values.
For $r\in[K]$ and $j\in\{0,\dots,r-1\}$, let
\begin{equation}\label{eq:mu-def}
  e_{r,j}
  :=
  \min_{I\subseteq[K]:\left|R_r\setminus I\right|=j}
  e_I;
\end{equation}
since $\left|R_r\setminus I_{\theta}\right|=j$ means that $R_r$ contains exactly $j$ true discoveries,
$e_{r,j}$ reflects the implausibility of there being exactly $j$ true discoveries among $R_r$.

\begin{lemma}\label{lem:frontier}
  For each $r\in[K]$,
  $\alpha_{R_r}=\alpha^*_r$ and $\beta_{R_r}=M_r/r$,
  where
  \begin{equation}\label{eq:M}
    M_r
    :=
    \max_{j\in\{0,\dots,r-1\}}
    \frac{r-j}{e_{r,j}}
  \end{equation}
  (with $c/0:=\infty$ for $c>0$).
  Therefore,
  \[
    \alpha^*_r
    =
    \begin{cases}
      M_r/r & \text{if $M_r\le r$}\\
      \infty & \text{otherwise}.
    \end{cases}
  \]
\end{lemma}

\begin{proof}
  The definitions immediately imply that $\beta_{R_r}=M_r/r$.
  It remains to prove that $\beta_R\ge\beta_{R_r}$ for every $R\subseteq[K]$ of size $r$.
  If $R\ne R_r$, there are $i_*<i^*$ such that $i_*\in R$ and $i^*\notin R$;
  set $R':=(R\setminus\{i_*\})\cup\{i^*\}$.
  It suffices to show that $\beta_{R'}\le\beta_R$,
  since repeating such replacements transforms $R$ into $R_r$.
  In other words, we need to check that for each $I\subseteq[K]$
  there exists $I'\subseteq[K]$
  such that $\FDP_I(R')/e_I\le\FDP_{I'}(R)/e_{I'}$.
  Let $I\subseteq[K]$.
  If $i^*\in I$ and $i_*\notin I$,
  we can set $I':=(I\setminus\{i^*\})\cup\{i_*\}$, so that
  $\left|R\cap I'\right|=\left|R'\cap I\right|$ and $e_{I'}\le e_I$.
  Otherwise, $\left|R\cap I\right|\ge\left|R'\cap I\right|$,
  and we can set $I':=I$.
\end{proof}

Let us fix $r\in[K]$.
For given sizes of $I\cap R_r$ and $I\setminus R_r$,
the arithmetic mean $e_I$ is smallest when these two sets consist of the indices of the smallest e-values
in $R_r$ and outside $R_r$, respectively.
Therefore, for each $j\in\{0,\dots,r-1\}$,
the minimum in \eqref{eq:mu-def} is attained at a set of the form $I=S_{r,j}\cup\{1,\dots,i\}$,
where
$
  S_{r,j}
  :=
  \{K-r+1,\dots,K-j\}
$
is $R_r$ without its $j$ largest elements
and $i\in\{0,\dots,K-r\}$.
Following \cite[Sect.~9]{Vovk/Wang:2023}, we set
\begin{align*}
  s_k &:= e_1+\dots+e_k, & k&=0,\dots,K,\\
  \sigma_j &:= e_{K-r+1}+\dots+e_{K-j} = s_{K-j}-s_{K-r}, & j&=0,\dots,r-1
\end{align*}
(the dependence of $\sigma_j$ on $r$ is not reflected in our notation),
so that
\begin{equation}\label{eq:mu}
  e_{r,j}
  =
  \min_{i\in\{0,\dots,K-r\}}
  \frac{s_i+\sigma_j}{i+r-j}.
\end{equation}
Computing $M_r$ directly using \eqref{eq:M} and \eqref{eq:mu} takes time $O(K^2)$,
which gives the overall time $O(K^3)$ for computing the frontier.
Algorithm~\ref{alg:frontier},
which is a minor modification of \cite[Algorithm~4]{Vovk/Wang:2023},
computes $M_r$ in time $O(K)$,
which gives the overall time $O(K^2)$.

\begin{algorithm}[bt]
  \caption{Computing the frontier}
  \label{alg:frontier}
  \begin{algorithmic}[1]
    \Require An increasing sequence of e-values $e_1\le\dots\le e_K$.
    \Ensure The frontier $\alpha^*_0,\dots,\alpha^*_K$.
    \State $s_0:=0$
    \For{$k=1,\dots,K$}
      \State $s_k:=s_{k-1}+e_k$
    \EndFor
    \State $\alpha^*_0:=0$
    \For{$r=1,\dots,K$}
      \State $k:=K-r$\label{ln:k_init}
      \State $M:=0$
      \For{$j=0,\dots,r-1$}
        \State $\sigma_j:=s_{K-j}-s_{K-r}$\label{ln:it_start}
        \State $\texttt{slope}:=\frac{s_k+\sigma_j}{k+r-j}$
        \For{$i=k-1,\dots,0$}
          \State $\texttt{new\_slope}:=\frac{s_i+\sigma_j}{i+r-j}$
          \If{$\texttt{new\_slope}>\texttt{slope}$} \textbf{break}\label{ln:break}
          \EndIf
          \State $k:=i$\label{ln:k}
          \State $\texttt{slope}:=\texttt{new\_slope}$
        \EndFor
        \State $M:=\max(M,(r-j)/\texttt{slope})$\label{ln:M} \Comment{$c/0:=\infty$}
      \EndFor
      \If{$M\le r$}
        \State $\alpha^*_r:=M/r$
      \Else
        \State $\alpha^*_r:=\infty$
      \EndIf
    \EndFor
  \end{algorithmic}
\end{algorithm}

\begin{figure}[bt]
  \begin{center}
    \includegraphics[width=0.55\textwidth]{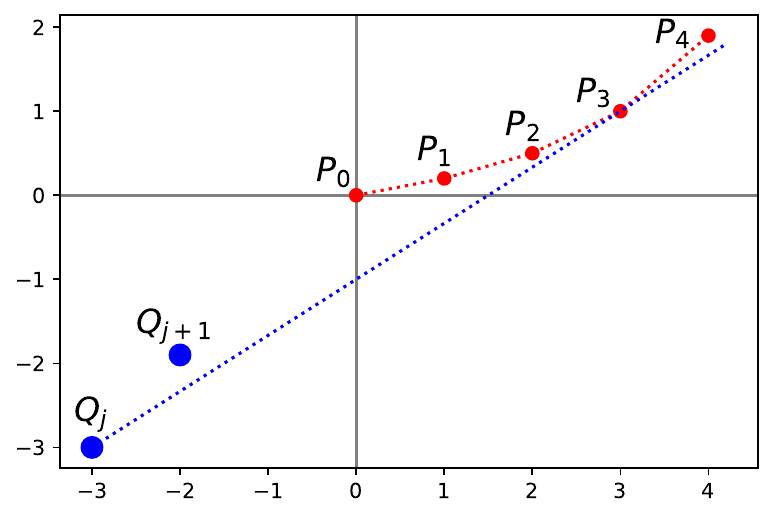}
  \end{center}
  \caption{Geometry behind Algorithm~\ref{alg:frontier}
    (reproduced from \cite[Fig.~8]{Vovk/Wang:2023})}
  \label{fig:algorithm}
\end{figure}

The geometry behind Algorithm~\ref{alg:frontier} is shown in Fig.~\ref{fig:algorithm}.
The coordinates of each of the points $P_k$, $k=0,\dots,K-r$, are $(k,s_k)$,
and the coordinates of the point $Q_j$, where $j\in\{0,\dots,r-1\}$, are $-(r-j,\sigma_j)$.
The slope of the straight line going from $Q_j$ to $P_i$ is
\[
  a_i
  :=
  \frac{s_i+\sigma_j}{i+r-j},
\]
and so $e_{r,j}$ is the smallest of these slopes over $i\in\{0,\dots,K-r\}$.
Since the sequence $e_1,\dots,e_{K-r}$ is increasing,
connecting the points $P_0,P_1,\dots$ in this order
(see the red line in Fig.~\ref{fig:algorithm})
gives us the graph of a convex function.
Therefore, $e_{r,j}$ is the slope of the straight line (shown in blue)
passing through $Q_j$ and touching the red line from below.
Let $P_{k_j}$ be the leftmost point at which the blue line touches the red line.
As $i$ increases from $0$ to $K-r$,
the slope $a_i$ strictly decreases until $i=k_j$ and does not decrease after that.

The command \textbf{break} in line~\ref{ln:break} of Algorithm~\ref{alg:frontier}
means leaving the innermost loop, as in Python or \textsf{R};
in this context, it is equivalent to ``go to line~\ref{ln:M}''.
The loops over $i$ are executed in the decreasing order of $i$
(and are not executed at all when $k=0$).
The variable $k$ in line~\ref{ln:k} is the index of the ``current vertex'' $P_k$;
we start from the rightmost $P_k$ in line~\ref{ln:k_init} and then keep moving left.
Figure~\ref{fig:algorithm} illustrates the execution of Algorithm~\ref{alg:frontier}
when $k=3$, so that the current vertex is $P_3$.
For each $j$, the loop over $i$ moves the current vertex to the left
as long as this does not make the line going from $Q_j$ to the current vertex steeper,
and line~\ref{ln:M} uses the slope of the resulting line,
which will be shown to be $e_{r,j}$.

\begin{theorem}
  Given the e-values \eqref{eq:sorted},
  Algorithm~\ref{alg:frontier} computes the frontier in time $O(K^2)$ and memory $O(K)$.
  If $r^*\in\{0,\dots,K\}$ maximizes $F(r,\alpha^*_r)$,
  the set $R_{r^*}$ is a solution of the optimization problem \eqref{eq:objective},
  and its level is $\alpha_{R_{r^*}}=\alpha^*_{r^*}$.
\end{theorem}

\begin{proof}
  Fix $r\in[K]$ and $j\in\{0,\dots,r-1\}$,
  and suppose that $k\ge k_j$ when the iteration of the loop over $j$ for this $j$
  starts in line~\ref{ln:it_start}.
  As long as $i\ge k_j$, we have $\texttt{new\_slope}\le\texttt{slope}$ in line~\ref{ln:break},
  and if $k_j>0$, the iteration with $i=k_j-1$ ends with \textbf{break}.
  Therefore, in line~\ref{ln:M} we have $k=k_j$ and $\texttt{slope}=e_{r,j}$.
  The assumption $k\ge k_j$ holds for $j=0$, since then $k=K-r$,
  and for $j>0$ it follows from $k_j\le k_{j-1}$, which we will now check.

  Let $j\in\{0,\dots,r-2\}$.
  First we check that the point $Q_{j+1}$ lies on or above the blue line passing through $Q_j$,
  as shown in Fig.~\ref{fig:algorithm}.
  The slope of the line going from $Q_j$ to $Q_{j+1}$ is $\sigma_j-\sigma_{j+1}=e_{K-j}\ge e_{K-r+1}$.
  If $k_j<K-r$, the slope of the blue line is at most $e_{k_j+1}\le e_{K-r}$,
  since the blue line touches the red line from below at $P_{k_j}$.
  If $k_j=K-r$, the slope of the blue line is the average of $e_1,\dots,e_{K-j}$,
  which does not exceed $e_{K-j}$.
  In both cases, the line going from $Q_j$ to $Q_{j+1}$ is at least as steep as the blue line,
  and so $Q_{j+1}$ lies on or above the blue line.
  Therefore, the line going from $Q_{j+1}$ to $P_{k_j}$ is not steeper than the blue line,
  and so
  $k_{j+1}\le k_j$.

  We can see that line~\ref{ln:M} replaces $M$ by $\max(M,(r-j)/e_{r,j})$,
  so that $M=M_r$ at the end of the loop over $j$,
  and $\alpha^*_r$ is computed correctly by Lemma~\ref{lem:frontier}.
  The variable $k$ never increases inside the loop over $j$,
  and each iteration of the loop over $i$ either decreases $k$ or is the last one for the given $j$;
  therefore, the loop over $j$ takes time $O(K)$.
  Computing $s_0,\dots,s_K$ takes time $O(K)$,
  and the memory used by the algorithm is $O(K)$.
  The last statement of the theorem follows from \eqref{eq:reduction} and Lemma~\ref{lem:frontier}.
\end{proof}

If the e-values are not sorted,
sorting them takes time $O(K\log K)$.
For a fixed $\alpha\in(0,1]$,
the largest $r$ such that $\alpha^*_r\le\alpha$ is the largest number of rejections
of the closed eBH procedure of \cite[(7)]{Xu/etal:arXiv2509} at level $\alpha$.
Algorithm~\ref{alg:frontier} finds it for all $\alpha$ at once in time $O(K^2)$,
whereas the algorithm in \cite[Sect.~10]{Xu/etal:arXiv2509}
takes time $O(K^2\log K)$ for one $\alpha$ in the worst case,
although, as its authors note, its running time is often closer to $K\log K$ in practice;
the algorithm mentioned in \cite[Sect.~9.6.1]{Ramdas/Wang:2025}
takes time $O(K^3)$.

\section{Conclusion}
\label{sec:conclusion}

This note discusses a procedure
that assigns a level to every rejection set
and so allows the researcher to choose the rejection set, together with its level,
after seeing the data.
The smallest level among the rejection sets of a given size
is attained by rejecting the base hypotheses with the largest e-values,
and the frontier can be computed in time $O(K^2)$.

A researcher who is interested in the false discovery rate
and is happy to fix the level $\alpha$ in advance
can simply use the closed eBH procedure at level $\alpha$;
see Sect.~\ref{sec:closed}
(although she may miss much larger rejection sets
available at a slightly larger level).
The relative eBH procedure does not have any parameters,
and its levels are valid in relative terms:
the expected ratio of the false discovery proportion to the reported level is at most $1$
(Proposition~\ref{prop:validity}).
It is an especially natural choice when the levels that the researcher might report
range widely
(overwhelming evidence can be reported with a very small level).

\subsection*{Acknowledgments}

It is a pleasure to thank Glenn Shafer, who turns 80 today,
for his encouragement and support over many years;
his ideas have shaped my thinking about the foundations of probability,
and about e-values in particular.

Claude Opus 5--5.5 and Fable 5.1--5.2 have been used
in exploring proof ideas and developing proofs,
which I reviewed carefully,
and in checking the note.
I take full responsibility for all claims and statements made here,
including mathematical statements and their proofs.

\appendix

\section{A large-deviation inequality}
\label{app:Bernstein}

The following lemma,
used in the proofs of Propositions~\ref{prop:closed-sequential} and \ref{prop:sequential},
follows immediately from a Bernstein-type inequality for supermartingales
as stated in \cite[Corollary~2.1]{Fan/etal:2012}.
In it, $\FFF_0\subseteq\FFF_1\subseteq\cdots$ is a filtration with a trivial $\FFF_0$,
such as the one introduced in Sect.~\ref{sec:closed};
notice that the random variables $Y_n$ are not required to be adapted to it.

\begin{lemma}\label{lem:Bernstein}
  Let $c\in(0,1]$, and let $E_n$, $c_n$, and $Y_n$, $n\in\N$, be random variables such that,
  for all $n\in\N$,
  $E_n$ is $\FFF_n$-measurable,
  $c_n$ is $\FFF_{n-1}$-measurable,
  $E_n\ge0$,
  $\E(E_n\mid\FFF_{n-1})\le1$,
  $c_n\in[0,c]$,
  and $Y_n\le\min(1,c_n(E_n-1))$.
  Then, for all $N\in\N$ and $\epsilon>0$,
  \begin{equation}\label{eq:Bernstein}
    \P
    \left(
      \frac1N
      \sum_{n=1}^N
      Y_n
      \ge
      \epsilon
    \right)
    \le
    \exp
    \left(
      -\frac{N\epsilon^2}{2c+2\epsilon/3}
    \right),
  \end{equation}
  and $\limsup_{N\to\infty}\frac1N\sum_{n=1}^N Y_n\le0$ a.s.
\end{lemma}

\begin{proof}
  Set $X_n:=\min(1,c_n(E_n-1))$.
  Then $Y_n\le X_n$, $-c\le X_n\le1$, and
  $\E(X_n\mid\FFF_{n-1})
  \le0$.
  Therefore, $X_1,X_2,\dots$ is a supermartingale difference sequence satisfying $-c\le X_n\le1$,
  and \cite[Corollary~2.1]{Fan/etal:2012} gives
  \[
    \P
    \left(
      \sum_{n=1}^N X_n\ge x
    \right)
    \le
    \exp
    \left(
      -\frac{2x^2}{4(N c+x/3)}
    \right)
  \]
  for all $x\ge0$.
  Since $\sum_{n=1}^N Y_n\le\sum_{n=1}^N X_n$,
  it remains to take $x:=N\epsilon$ to obtain \eqref{eq:Bernstein}.
  Finally, \eqref{eq:Bernstein} and the Borel--Cantelli lemma show that,
  for each $\epsilon>0$, almost surely $\frac1N\sum_{n=1}^N Y_n<\epsilon$ from some $N$ on,
  which implies the last statement of the lemma.
\end{proof}
\end{document}